\documentclass[12pt]{amsart}

\usepackage{graphicx,tikz}
\usepackage{latexsym}
\usepackage{amsmath}
\usepackage{amssymb}
\usepackage{array}
\usepackage{arydshln}

\usepackage{cite}
\usepackage{stmaryrd}
\usepackage{enumerate}

\usepackage{hyperref}

\usepackage{color}
\usepackage{lineno}
\usepackage{graphicx}
\usepackage{ae}
\usepackage{amsmath}
\usepackage{amssymb}
\usepackage{latexsym}
\usepackage{url}
\usepackage{epsfig}
\usepackage{mathrsfs}
\usepackage{amsfonts}
\usepackage{amsthm}
\usepackage{bm}

\newcommand{\Cay}{\mathrm{Cay}}

\newcommand{\Sym}{\mathrm{Sym}}
\newcommand{\Orb}{\mathrm{Orb}}
\newcommand{\fix}{\mathrm{fix}}
\newcommand{\cn}{\mathsf{c}}
\newcommand{\on}{\mathsf{o}}
\newcommand{\Aut}{\mathrm{Aut}}

\newcommand{\cO}{\mathcal{O}}
\newcommand{\D}{\mathrm{D}}
\newcommand{\ZZ}{\mathbb{Z}}
\def\mod{{\rm mod\,}}

\newtheorem{theorem}{Theorem}[section]

\newtheorem{lemma}[theorem]{Lemma}

\newtheorem{cor}[theorem]{Corollary}

\theoremstyle{definition}

\numberwithin{equation}{section} 
\allowdisplaybreaks    

\def\proof{\par\noindent{\bf Proof.~}}

\def\qed{\hfill$\Box$\vspace{12pt}}

\long\def\delete#1{}

\usepackage{xcolor}
\usepackage[normalem]{ulem}

\usepackage{stfloats}

\begin{document}
\title{Enumerating Cayley digraphs on dihedral groups}
\thanks{Supported by the National Natural Science Foundation of China (12471328, 12331013) and the Fundamental Research Funds for the Central Universities.}

\author[Lu]{Zai Ping Lu}
\address{
Zaiping Lu\\
Center for Combinatorics, LPMC,
Nankai University,
Tianjin 300071, China
}
\email{lu@nankai.edu.cn}
\author[Xie]{Jia Yin Xie}
\address{
Jiayin Xie\\
Center for Combinatorics, LPMC,
Nankai University,
Tianjin 300071, China
}
\email{xiejiayin2022@163.com}
\author[Xie]{Jin-Hua Xie}
\address{Jin-Hua Xie\\
Center for Combinatorics, LPMC,
Nankai University,
Tianjin 300071, China}
\email{jinhuaxie@nankai.edu.cn}

\date{}

\openup 0.5\jot
\maketitle
\begin{abstract}
This paper investigates the enumeration of Cayley digraphs, focusing on counting   Cayley digraphs on dihedral groups up to CI-isomorphism. By leveraging the Cauchy-Frobenius Lemma and properties of automorphisms, we derive an explicit formula  for the number of non-isomorphic Cayley digraphs on dihedral groups with DCI-property, particularly for the group $\D_{6p}$
with $p>3$  a prime. The enumeration involves detailed analysis of cycle numbers of automorphisms and their actions on the group elements, culminating in a precise count of non-isomorphic digraphs.

\vskip 10pt

\end{abstract}

\vskip 20pt

\section{Introduction}

Cayley (di)graphs serve as a fundamental tool in algebraic graph theory, establishing a deep connection between group theory and graph theory. Extensive research has been devoted to their structural properties, including isomorphism classification, symmetry characteristics, and automorphism groups, leading to a wealth of theoretical results. However, the enumeration of non-isomorphic Cayley (di)graphs, a problem of significant interest yet considerable difficulty, remains comparatively underdeveloped, presenting an important challenge in the field.

A Cayley digraph $\Cay(R,S)$ is constructed from a group
$R$ and a subset $S\subseteq R\setminus\{1\}$ (called a connection set),  where the vertices represent the elements of $R$, and directed edges correspond to left multiplication by elements of $S$. If $S$ is inverse-closed (i.e., $S=S^{-1}:=\left\{s^{-1}\mid s\in S\right\}$), then $\Cay(R,S)$ can be regarded as an undirected graph, in which case it is called  a  Cayley graph on $R$. Clearly, if
$a$ is an automorphism  of the group $R$ then $\Cay(R,S)$ and $\Cay(R,S^a)$ are isomorphic.
Two Cayley digraphs $\Cay(R,S)$ and $\Cay(R,T)$ are said to be CI-isomorphic if there exists an automorphism $a$ of the group $R$ such that $T=S^a$. In this case,  the connection sets $S$ and $T$ are called CI-equivalent. A group $R$ is called a DCI-group (or CI-group) if every pair of isomorphic Cayley digraphs (or Cayley graphs) of $R$ must be CI-isomorphic.

In 1967, Turner \cite{Turner1967}   showed  that every vertex-transitive graph of prime order
$p$  is necessarily a Cayley graph on the cyclic group $\ZZ_p$, and two such Cayley graphs are isomorphic if and only if their connection  sets are CI-equivalent. Leveraging this characterization, Turner applied P{\'o}lya enumeration theorem to derive the precise count of non-isomorphic vertex-transitive graphs of order $p$.
Alspach and Mishna \cite{AlspachM2002} extended Turner's enumeration result to the case of Cayley digraphs and Cayley graphs on cyclic groups. By applying  P{\'o}lya enumeration theorem, they developed a systematic counting method for non-isomorphic Cayley digraphs on DCI-groups and non-isomorphic Cayley graphs on CI-groups. In particular, based on the complete classification of cyclic DCI-groups and CI-groups \cite{Muzychuk1995,Muzychuk1997,Muzychuk1997-1}, they precisely determined the numbers of non-isomorphic Cayley digraphs and Cayley graphs on the cyclic group $\ZZ_n$, where $n$ is either square-free or twice a square-free integer.
More recently,
given that the dihedral group $\D_{2p}$  is a DCI-group \cite{Babai1977},
Huang and Huang \cite{HuangH2019} extended these results to
$\D_{2p}$,  where $p$ is a prime, providing a detailed count of non-isomorphic Cayley digraphs and graphs  on $\D_{2p}$. For more related results, please refer to \cite{Chia-Ong,HuangHL2017,LiSim2001,LiskovetsP2000,Mamut}.

In this paper, we focus on the dihedral group  $\D_{2n}$ of order $2n$, which models the symmetries of a regular
$n$-gon, is a particularly interesting case due to its rich structure and widespread applications. A key challenge lies in determining when two Cayley digraphs on
$\D_{2n}$ are isomorphic. As for this,  we currently have no good ideas. Here, we only consider those dihedral groups that are DCI-groups, which reduces to understanding the action of the automorphism group on the power set of the group's non-identity elements.
Our primary goal is to enumerate the non-isomorphic Cayley digraphs on $\D_{6p}$, where $p>3$ is a prime.

The paper is organized as follows: Section \ref{sect2} reviews the necessary background on group actions and the Cauchy-Frobenius Lemma. Section \ref{sect=autoD2n} details the automorphisms of dihedral groups and their properties. Section \ref{sect=cyclenum} derives the cycle numbers of these automorphisms, which are crucial for the enumeration. Finally, Section \ref{sect=D6p} applies these results to the specific case of
$\D_{6p}$, providing an explicit formula  for the number of non-isomorphic Cayley digraphs.

\vskip 20pt

\section{Enumeration of orbits}\label{sect2}

Let $G$ be a finite group acting on a  finite nonempty set   $\Omega$, that is,
there exists a map
\[\Omega\times G\rightarrow \Omega,\,\,(\alpha, x)\mapsto \alpha^x\] such that
\begin{itemize}
\item[(i)] $\alpha=\alpha^1$ for all $\alpha \in \Omega$, where $1$ is the identity element of $G$,
\item[(ii)] $\alpha^{xy}=(\alpha^x)^y$ for all $x,y\in G$ and all $\alpha\in \Omega$.
\end{itemize}
 For a point  $\alpha\in \Omega$, put
$\Orb_\Omega(\alpha,G): =\{\alpha^x\mid x\in G\}$,
called the {\em orbit} of $\alpha$ under $G$, or a {\em $G$-orbit}. Clearly, $\Orb_\Omega(\beta,G)=\Orb_\Omega(\alpha,G)$ for any $\beta\in \Orb_\Omega(\alpha,G)$.
Then $\Omega$ is partitioned into distinct $G$-orbits, say \[\Omega=\biguplus_{i=1}^m\Orb_\Omega(\alpha_i,G),\]
where $m$ is called the {\em orbit number} of $G$ (on $\Omega$), denoted by $\on(\Omega,G)$.
Denote $G_\alpha$  the stabilizer  of $\alpha$ in $G$, that is,
$G_\alpha:=\{x\in G\mid \alpha^x=\alpha\}$,
which is a subgroup of $G$.
Then, for  $\beta\in \Orb_\Omega(\alpha,G)$, we have \[|\Orb_\Omega(\alpha,G)||G_\alpha|=|G|=|\Orb_\Omega(\alpha,G)||G_\beta|,\] in particular, $|G_\alpha|=|G_\beta|$. In fact, $G_\beta$ is a conjugation of $G_\alpha$ in $G$, that is, if $\beta=\alpha^x$ then $G_\beta=x^{-1}G_\alpha x$.
Enumerating those pairs $(\gamma,x)\in \Omega\times G$ with  $\gamma^x=\gamma$,
\[\sum_{x \in G} |\fix_\Omega(x)|=\sum_{\gamma\in \Omega}|G_\gamma|=\sum_{i=1}^m\sum_{\gamma\in \Orb_\Omega(\alpha_i,G)}|G_\gamma|=\sum_{i=1}^m|\Orb_\Omega(\alpha_i,G)||G_{\alpha_i}|=m|G|,\]
where $\fix_\Omega(x):=\{\gamma\in \Omega\mid \gamma^x=\gamma\}$.
Then we have the famous Cauchy-Frobenius Lemma, refer to \cite[p. 24, Theorem 1.7A ]{Dixon}.

\begin{lemma}[Cauchy-Frobenius Lemma]\label{enlemma1}
Let $ G $ be a finite group acting on a finite nonempty set $ \Omega $. Then the number  of $G$-orbits is
\[
\on(\Omega,G) = \frac{1}{|G|}\sum_{x \in G} |\fix_\Omega(x)|.
\]
\end{lemma}

By the above (i) and (ii),  it is easy to see that every element $x\in G$ induces a permutation ${\bar x}$ of $\Omega$, and the mapping $x\mapsto {\bar x}$ is a  homomorphism from the group $G$ to the symmetric group $\Sym(\Omega)$ on $\Omega$. We denote $\cn(\Omega,x)$ the number of cycles (including cycles of length $1$) appearing in a factorization of ${\bar x}$ into disjoint cycles, called the {\em cycle number} of $x$ (on $\Omega$). Note that
$\cn(\Omega,x)$ is just the number $\on(\Omega,\langle x\rangle)$ of   orbits on $\Omega$ of the cyclic group $\langle x\rangle$.
Then, by  Cauchy-Frobenius Lemma, we have the following simple fact.

\begin{lemma}\label{enlemma3}
Let $G$ be a finite group acting on a finite nonempty set $\Omega$, and let $x\in G$ have order $|x|$. Then  cycle number of $x$ is equal to
 \[        \cn(\Omega,x)  = \on(\Omega,\langle x\rangle)  = \frac{1}{|x|} \sum_{y \in \langle x \rangle} |\fix_\Omega(y)|.
    \]
\end{lemma}

Denote $2^\Omega$ the power set of $\Omega$. Then $G$ acts on $2^\Omega$, by
\[S^x:=\{\alpha^x\mid \alpha \in S\},\quad x\in G,\,\, S\in 2^\Omega.\]
Clearly, a subset $S\subseteq\Omega$ is fixed by an element $x\in G$ if and only of
$S$ is the union of some $\langle x\rangle$-orbits on $\Omega$. This implies that
\[|\fix_{2^\Omega}(x)|=2^{\cn(\Omega,x)},\quad \mbox{for all }x\in G.\]
Using the Cauchy-Frobenius Lemma,
one may formulate $\on(2^\Omega, G)$ as follows.

\begin{lemma}\label{enlemma2}
Let $G$ be a finite group acting on a finite nonempty set $\Omega$. Then the number of $G$-orbits on $2^\Omega$ is equal to
\[
\on(2^\Omega, G)=\frac{1}{|G|} \sum_{x \in G} 2^{\cn(\Omega,x)}.
\]
\end{lemma}

In the following, $R\ne 1$ is assumed to be a finite group,  and put $R^\#: =R\setminus \{1\}$.
Let $A:=\Aut(R)$ be the automorphism group of $R$.   Clearly, $A$ is a subgroup of $\Sym(R^\#)$.
Denote $\cO$ the set of $A$-orbits on $2^{R^\#}$.
Two subsets $S,T\subseteq R^\#$ are called {\em CI-equivalent} if they are contained in a same member of $\cO$, that is, there exists $a\in A$ such that $T=S^a$.
If $S$ and $T$ are CI-equivalent, then every element $a\in A$ with $S^a=T$ gives an isomorphism from
  $\Cay(R,S)$ to $\Cay(R,T)$; in this case, we say that  $\Cay(R,S)$ and $\Cay(R,T)$ are {\em CI-isomorphic}. Thus, up to CI-isomorphism,
  there are exactly $\on(2^{R^\#}, A)$ many Cayley digraphs on $R$. By Lemma \ref{enlemma2},
    \begin{equation}\label{eq-o}
    \on(2^{R^\#}, A)=\frac{1}{|A|} \sum_{a \in A} 2^{\cn(R^\#,a)}.
  \end{equation}

In general, an isomorphism between two Cayley digraphs   $\Cay(R,S)$ and $\Cay(R,T)$ is not necessary a  CI-isomorphism. Thus, $ \on(2^{R^\#}, A)$ is greater than or equal to the number of Cayley digraphs  on $R$ up to isomorphism, where the equality holds if and only if subsets of $R^\#$ in distinct members of $\cO$ give rise to  non-isomorphic Cayley digraphs.

Recall that a subset $S\subseteq R^\#$ is called a {\em CI-subset} of $R$ if  $\Cay(R,T)\cong \Cay(R,S)$ for any $T\subseteq R^\#$, then $S$ and $T$ are CI-equivalent, while the Cayley digraph $\Cay(R,S)$ is called a {\em CI-digraph}, and called a {\em CI-graph} if further $S=S^{-1}$.
The group $R$ is called a {\em DCI-group} if every subset $S\subseteq R^\#$ is a CI-subset, and called a {\em CI-group} if every subset $S\subseteq R^\#$ with $S=S^{-1}$ is a CI-subset.

By the argument above, we have the following simple fact, see also \cite[Theorem 1.6]{AlspachM2002}.

\begin{lemma}\label{enlemma2-1}
Let $R\ne 1$   be a finite group,    $R^\#: =R\setminus \{1\}$ and $A:=\Aut(R)$.
If $R$ is a DCI-group then the number of non-isomorphic Cayley digraphs on $R$ is given as in \eqref{eq-o}.
\end{lemma}

\vskip 20pt

\section{Automorphisms of   dihedral groups}\label{sect=autoD2n}
In this section, $n$ is always assumed to be a positive integer.

Put
$\D:=\langle u,v\mid v^2=1, vuv=u^{-1}\rangle$,
 the infinite dihedral group.
Then  the  finite dihedral group $\D_{2n}$ of order $2n$ is a  homomorphic image of $\D$
via the natural homomorphism form $\D$ to $\D/\langle u^n\rangle$. We denote $u_n$ and $v_n$ the images of $u$ and $v$ respectively, and write
\[\D_{2n}:=\langle u_n,v_n\rangle.\]
Then $\langle u_n\rangle$ is the unique subgroup of order $n$ in $\D_{2n}$, and every element of $\D_{2n}$
is contained in either $\langle u_n\rangle$ or the coset $\langle u_n\rangle v_n$.

Denote $\ZZ$ the additive group of integers.
For each pair $r,t\in \ZZ$, define a mapping $a_{r,t}: \D\rightarrow \D$ by
\[ (u^i)^{a_{r,t}}:= u^{ri},\, \, (u^jv)^{a_{r,t}}:= u^{rj+t}v.\]
Then $a_{r,t}$ is a homomorphism, and
\begin{equation}\label{eq-anrs}
a_{n,r,t}: \D_{2n}\rightarrow \D_{2n},\,\,  u_n^i\mapsto u_n^{ri},\, \, u_n^jv_n\mapsto u_n^{rj+t}v_n
\end{equation}
is also  a homomorphism. Moreover,
\begin{equation}\label{eq-a_{r,t}^s}
  \begin{split}
  &a_{r,t}a_{r',t'}=a_{rr', r't+t'},\\
  &a_{n,r,t}a_{n,r',t'}=a_{n,rr', r't+t'},\\
  &a_{n,r,t}=a_{n,r',t'}\Leftrightarrow r\equiv r'\,\,(\mod n), t\equiv t'\,(\mod n).
  \end{split}
\end{equation}
It is easy to see that every element in $\Aut(\D_{2n})$ has the form of $a_{n,r,t}$ with $\gcd(n,r)=1$. Thus
\begin{equation}\label{eq-AutD_2n}
 \Aut(\D_{2n})=\{a_{n,r,t}\mid r, t\in \ZZ, \gcd(n,r)=1\}.
\end{equation}
In view of \eqref{eq-a_{r,t}^s}, we have
\[|\Aut(\D_{2n})|=n\phi(n),\]
where $\phi$ is the Euler's phi function.
In addition, the subscripts $r$ and $t$ of an automorphism  $a_{n,r,t}\in \Aut(\D_{2n})$  are always taken modulo   $n$.

In the following, we consider the orders of elements in $\Aut(\D_{2n})$.
For $r\in \ZZ$ with $\gcd(r,n)=1$, we denote
 $|r|_n$ the order of $r$ modulo $n$. Clearly, if $d\geqslant 1$ is a divisor of
$n$ then $|r|_d$ is a divisor of $|r|_n$.

\begin{lemma}\label{r_n-d}
Let $n$, $m$ and $d$ be positive integers with $n=dm$ and $\gcd(d,m)=1$. If   $r\in \ZZ$ with $\gcd(r,n)=1$ then
$|r|_n=\frac{|r|_d|r|_m}{\gcd(|r|_d,|r|_m)}$.
\end{lemma}
\proof
Let $r\in \ZZ$ with $\gcd(r,n)=1$. Clearly, $|r|_n$ is a common multiple of $|r|_d$ and $|r|_m$. Let $k=\frac{|r|_d|r|_m}{\gcd(|r|_d,|r|_m)}$. Then $r^k-1$ has divisors $r^{|r|_d}-1$ and $r^{|r|_m}-1$,
and so $r^k-1\equiv 0\,(\mod d)$ and  $r^k-1\equiv 0\,(\mod m)$. Since $n=dm$ and $\gcd(d,m)=1$, we have
$r^k-1\equiv 0\,(\mod n)$. This implies that $|r|_n$ is a divisor of $k$. Then $|r|_n=k=\frac{|r|_d|r|_m}{\gcd(|r|_d,|r|_m)}$, as desired.
\qed

Define a polynomial of degree $n-1$
\[S_n(x):=1+x+\cdots+x^{n-1},\]
and put $S_0(x):=0$. For $r,\,t \in \ZZ$ with $\gcd(r,n)=1$, it is easily seen from \eqref{eq-a_{r,t}^s} that
\begin{equation}\label{eq-a_{r,t}^s-1}
 a_{n,r,t}^s=a_{n,r^s,tS_{s}(r)} \mbox{ for every integer } s\geqslant 1.
\end{equation}
Next we define
 \[\kappa(n,r,t):=\frac{n|r|_n}{\gcd\left(n, tS_{|r|_n}(r)\right)}.\]
We will show in Lemma~\ref{lem-order-a-rt} that $|a_{n,r,t}|=\kappa(n,r,t)$. To prove this, we need the following two lemmas.

\begin{lemma}\label{lem-kappa-eq-1}
Let $t,\,t',\,r\in \ZZ$ with $\gcd(r,n)=1$.
\begin{enumerate}[\rm(i)]
  \item If $\gcd(t,n)=\gcd(t',n)$ then $\kappa(n,r,t)=\kappa(n,r,t')$.
  \item If $t\mid t'$ then $\kappa(n,r,t')\mid \kappa(n,r,t)$, and $\kappa(n,r,t)=\kappa(n,r,t')\gcd\left(\frac{ \kappa(n,r,t)}{|r|_n},\frac{t'}{t}\right)$.
  \item If $d\geqslant 1$ is a divisor of
$n$ then, for all integers $\ell, k\geqslant 0$, $S_{\ell k|r|_d}(r)\equiv \ell S_{ k|r|_d}(r)\,(\mod d)$.
\end{enumerate}
\end{lemma}
\proof
All the statements can be verified easily.\qed

\begin{lemma}\label{lem-kappa}
Let $n$, $m$ and $d$ be positive integers with $n=dm$ and $\gcd(d,m)=1$. Let   $r,\,t\in \ZZ$ with $\gcd(r,n)=1$.  Then
\[
    \kappa(n,r,t)=\frac{|r|_n\kappa(d,r,t)\kappa(m,r,t)}{\gcd\left(\kappa(d,r,t),|r|_n\right)\gcd\left(\kappa(m,r,t),|r|_n\right)}.
\]
In particular, $\kappa(d,r,t)$ is a divisor of $\kappa(n,r,t)$.
 \end{lemma}
 \proof
Let $k=\frac{|r|_m}{\gcd(|r|_d,|r|_m)}$ and  $\ell=\frac{|r|_d}{\gcd(|r|_d,|r|_m)}$. Since $n=dm$ and $\gcd(d,m)=1$, by Lemma \ref{r_n-d}, $|r|_n=|r|_dk=|r|_m\ell$. Then, by Lemma~\ref{lem-kappa-eq-1}(iii),
 \[
 tS_{|r|_n}(r)\equiv ktS_{|r|_d}(r)\,(\mod d),\,\,\, tS_{|r|_n}(r)\equiv \ell tS_{|r|_m}(r)\,(\mod m).
 \]
  This implies that
 \[\gcd\left(d,tS_{|r|_n}(r)\right)=\gcd\left(d, ktS_{|r|_d}(r)\right),\,\, \gcd\left(m,tS_{|r|_n}(r)\right)=\gcd\left(m, \ell tS_{|r|_m}(r)\right).\]
 Hence, it follows from $\gcd(d,m)=1$ that
 \begin{align*}
      \gcd\left(n,tS_{|r|_n}(r)\right)&=\gcd\left(d,tS_{|r|_n}(r)\right)\gcd\left(m,tS_{|r|_n}(r)\right)\\
      &=\gcd\left(d,ktS_{|r|_d}(r)\right)\gcd\left(m,\ell tS_{|r|_m}(r)\right).
 \end{align*}
 Then we derive from Lemma~\ref{lem-kappa-eq-1}(ii) that
 \[\begin{split}
 \gcd(|r|_d,|r|_m)\kappa(n,r,t)&=\gcd(|r|_d,|r|_m)\frac{n|r|_n}{\gcd\left( n,tS_{|r|_n}(r)\right) }\\
 &=\frac{d|r|_dm|r|_m}{\gcd\left(d, ktS_{|r|_d}(r)\right)\gcd\left(m, \ell tS_{|r|_m}(r)\right)}\\
 &=\kappa\left( d,r,tk\right)
   \kappa\left( m,r,t\ell\right)\\
   &=\frac{\kappa(d,r,t)}{\gcd\left(\frac{\kappa(d,r,t)}{|r|_d},k\right)}
   \frac{\kappa(m,r,t)}{\gcd\left(\frac{\kappa(m,r,t)}{|r|_m},\ell\right)}\\
   &=\frac{\kappa(d,r,t)\kappa(m,r,t)|r|_d|r|_m}{\gcd\left(\kappa(d,r,t),|r|_n\right)\gcd\left(\kappa(m,r,t),|r|_n\right)}.
 \end{split}\]
Therefore, we have
 \[\begin{split}
 \kappa(n,r,t)&=\frac{\kappa(d,r,t)\kappa(m,r,t)}{\gcd\left(\kappa(d,r,t),|r|_n\right)\gcd\left(\kappa(m,r,t),|r|_n\right)}
 \frac{|r|_d|r|_m}{\gcd(|r|_d,|r|_m)}\\
 &=\frac{|r|_n\kappa(d,r,t)\kappa(m,r,t)}{\gcd\left(\kappa(d,r,t),|r|_n\right)\gcd\left(\kappa(m,r,t),|r|_n\right)},
 \end{split}
 \]
 as required.\qed

\begin{lemma}\label{lem-order-a-rt} Let $a_{n,r,t}\in  \Aut(\D_{2n})$. Then $|r|_n$ is a divisor of $|a_{n,r,t}|$, and \[|a_{n,r,t}|=\frac{n|r|_n}{\gcd\left(n,tS_{|r|_n}(r)\right)}=\kappa(n,r,t).\]
In particular, $tS_{\kappa(n,r,t)}(r)\equiv 0\,(\mod n)$,
and if $\gcd(r-1,n)=1$ then $|a_{n,r,t}|=|r|_n$.
\end{lemma}
\proof
It is clear from \eqref{eq-a_{r,t}^s-1} that if $a_{n,r,t}^s=a_{n,1,0}$ then   $r^s\equiv 1\,(\mod n)$. Noting that $a_{n,1,0}$ is the identity element of $\Aut(\D_{2n})$, it follows that $|r|_n$ is a divisor of $|a_{n,r,t}|$.
Now let $m:=|a_{n,r,t}|$, and write $m=k|r|_n$. Then
\[
a_{n,1,0}=a_{n,r,t}^m=a_{n,r,t}^{k|r|_n}=a_{n,1,tS_{k|r|_n}(r)},
\]
which combined with Lemma~\ref{lem-kappa-eq-1}(iii) implies that
\[
a_{n,1,0}=a_{n,1,ktS_{|r|_n}(r)}.
\]
Thus, $k$ is the smallest positive integer such that $ktS_{|r|_n}(r)\equiv 0\,(\mod n)$, and so $k=\frac{n}{\gcd\left(n,tS_{|r|_n}(r)\right)}$. Then
\[
|a_{n,r,t}|=m=k|r|_n=\frac{n|r|_n}{\gcd\left(n,tS_{|r|_n}(r)\right)}=\kappa(n,r,t).
\]
Note that $n\mid (r^{|r|_n}-1)$ and $r^{|r|_n}-1=(r-1)S_{|r|_n}$. If $\gcd(r-1,n)=1$ then $S_{|r|_n}(r)\equiv 0\,(\mod n)$, yielding $k=1$. Thus the lemma follows.
\qed

\begin{cor}\label{lem-kappa-eq=2}
Let $t,\,y ,\,r, \, x\in \ZZ$ with $\gcd(r,n)=1=\gcd(x,n)$. Then $\kappa(n,r,t)=\kappa\left(n,r,(1-r)y+xt\right)$ and
$\gcd\left(n, tS_{|r|_n}(r)\right)=\gcd\left(n,((1-r)y+xt)S_{|r|_n}(r)\right)$.
 \end{cor}
 \proof
 Let $k=\kappa(n,x,y)$. Then by \eqref{eq-a_{r,t}^s-1} and Lemma~\ref{lem-order-a-rt}, we have
 \[
a_{n,x,y}^{-1}=a_{n,x,y}^{k-1}=a_{n, x^{k-1}, yS_{k-1}(x)}.
\]
This together with \eqref{eq-a_{r,t}^s} gives that
\[
a_{n,x,y}^{-1}a_{n,r,t}a_{n,x,y}=a_{n,x^{k-1},yS_{k-1}(x)}a_{n,r,t}a_{n,x,y}=a_{n,r,rxyS_{k-1}(x)+xt+y}.
\]
Noting that $xS_{k-1}(x)=x+\cdots+x^k=S_{k}(x)-1$, we obtain from Lemma~\ref{lem-order-a-rt} that
\[
rxyS_{k-1}(x)=ryS_k(x)-ry \equiv-ry\,(\bmod\,n).
\]
Hence, we drive from \eqref{eq-a_{r,t}^s} that
\[
a_{n,x,y}^{-1}a_{n,r,t}a_{n,x,y}=a_{n,r,rxyS_{k-1}(x)+xt+y}=a_{n,r,(1-r)y+xt}.
\]
Since $a_{n,r,t}$ and $a_{n,x,y}^{-1}a_{n,r,t}a_{n,x,y}$ have the same order, we conclude from Lemma \ref{lem-order-a-rt} that
\[\kappa(n,r,t)=|a_{n,r,t}|=|a_{n,r,(1-r)y+xt}|=\kappa(n,r,(1-r)y+xt),\]
as required.
\qed

\vskip 20pt

\section{Cycle numbers  of   automorphisms}\label{sect=cyclenum}
Let $n>1$ be an integer. Recall that
\[\D_{2n}:=\langle u_n,v_n\mid u_n^n=1=v_n^2, v_nu_nv_n=u_n^{-1}\rangle.\]
For each divisor $d\geqslant 1$ of $n$, denote $U_n(d)$ the set of elements with order $d$ in $\langle u_n\rangle$.
Then
\[U_n(d)=\{u_n^{\frac{ni}{d}}\mid \gcd(i,d)=1\}.\]
Clearly, $|U_n(d)|=\phi(d)$.
Put \[V_n:=\{v_n,u_nv_n,u_n^2v_n,\ldots, u_n^{n-1}v_n\}.\]
We has a partition of $\D_{2n}^\#$ shown as follows:
\[\D_{2n}^\#=V_n\uplus \biguplus_{d>1\text{ and }d \mid n }U_n(d).\]
Put
\[
U_n:=\biguplus_{d>1\text{ and }d \mid n}U_n(d).
\]
Since $V_n$ and $U_n(d)$ are invariant under $\Aut(\D_{2n})$ for each divisor $d\geq 1$ of $n$, we derive from the definition of $\cn(\Omega,x)$ before Lemma~\ref{enlemma3} that, for every $a_{n,r,t}\in \Aut(\D_{2n})$,
\begin{equation}\label{eq-two-part-c(Dn,art)}
 \cn(\D_{2n}^\#, a_{n,r,t})=\cn(V_n, a_{n,r,t})+\cn(U_n, a_{n,r,t})=\cn(V_n, a_{n,r,t})+\sum_{d>1\text{ and }d \mid n}\cn(U_n(d), a_{n,r,t}).
\end{equation}

\begin{lemma}\label{lem-Un(d)}
Let   $a_{n,r,t}\in \Aut(\D_{2n})$. Then
\[\cn(U_n(d), a_{n,r,t})=\cn(U_n(d), a_{n,r,0})= \frac{\phi(d)}{|r|_d}\]
for every divisor $d>1$ of $n$, and so
\[
\cn(U_n, a_{n,r,t})=\cn(U_n, a_{n,r,0})=\sum_{d>1\text{ and }d \mid n} \frac{\phi(d)}{|r|_d}.
\]
\end{lemma}
\proof
Clearly, $(u^i)^{a_{n,r,t}}=(u^i)^{a_{n,r,0}}$ for all $u_n^i\in U_n(d)$, and so $\cn(U_n(d), a_{n,r,t})=\cn(U_n(d), a_{n,r,0})$.
For $1\leqslant s\leqslant |a_{n,r,0}|=|r|_n$, we have $a_{n,r,0}^s=a_{n,r^s,0}$.
Noting that $U_n(d)=\{u_n^{\frac{ni}{d}}\mid \gcd(i,d)=1\}$, we have
\[
\begin{split}
u_n^{\frac{ni}{d}}\in \fix_{U_n(d)}(a_{n,r,0}^s)&\Leftrightarrow \frac{n}{d}i(r^s-1)\equiv 0\,(\mod n)\\
&\Leftrightarrow r^s-1\equiv 0\,(\mod d)\\
&\Leftrightarrow |r|_d\mid \gcd(s,|r|_n)\\
&\Leftrightarrow s=k|r|_d \mbox{ with } 1\leqslant k\leqslant \frac{|r|_n}{|r|_d}.
\end{split}
\]
Thus
\[
|\fix_{U_n(d)}(a_{n,r,0}^s)|=\begin{cases}
|U_n(d)|=\phi(d) &\text{ if }s=k|r|_d\text{ with }\,1\leqslant k\leqslant \frac{|r|_n}{|r|_d}\\
0 &\text{ otherwise }.\\
\end{cases}
\]
Clearly, $|a_{n,r,0}|=|r|_n$. By Lemma \ref{enlemma3}, we have

\[
\cn(U_n(d), a_{n,r,0}) = \frac{\phi(d)}{|r|_n} \left|\left\{k\,\Big |\,1\leqslant k\leqslant \frac{|r|_n}{|r|_d}\right\}\right|=\frac{\phi(d)}{|r|_d},
 \]
    and the lemma follows.
\qed

 We next formulate  the number $\cn(V_n, a_{n,r,t})$. Clearly, $V_n$ generates the dihedral group $\D_{2n}$.
 Thus, if some element in  $\Aut(\D_{2n})$ fixes  $V_n$ pointwise, then it is the identity automorphism $a_{n,1,0}$.  This implies that the order  $|a_{n,r,t}|$ of  $a_{n,r,t}\in \Aut(\D_{2n})$ is the least common multiple of the lengthes of cycles appearing in a cycle factorization of $a_{n,r,t}$ (on $V_n$).
 Thus, if $\gcd(s,|a_{n,r,t}|)=1$ then $|\fix_{V_n}(a_{n,r,t})|=|\fix_{V_n}(a_{n,r,t}^s)|$.
 In particular, two elements in $\langle a_{n,r,t}\rangle$ with the same order fix  the same number of
 points in $V_n$. We deduce
 \[
 \sum_{a \in \langle a_{n,r,t} \rangle} |\fix_{V_n}(a)|=\sum_{s\mid \kappa(n,r,t)}\phi\left( \frac{  \kappa(n,r,t)}{s}\right)  |\fix_{V_n}( a_{n,r,t}^s)|.
 \]
  By Lemma \ref{enlemma3}, we have

 \begin{equation}\label{eq-cV_n}
 \cn(V_n, a_{n,r,t})= \frac{1}{ \kappa(n,r,t)}\sum_{s\mid  \kappa(n,r,t)}\phi\left( \frac{ \kappa(n,r,t)}{s}\right)  |\fix_{V_n}( a_{n,r,t}^s)|.
 \end{equation}

 For $u_n^iv_n\in V_n$, we have $(u_n^iv_n)^{a_{n,r,t}}=u_n^{ri+t}v_n$, and so
\[u_n^iv_v\in \fix_{V_n}(a_{n,r,t})\Leftrightarrow u_n^{ri+t}v_n=u_n^iv_n\Leftrightarrow (r-1)i\equiv -t\,(\mod n).\]
By a result on the solutions of linear congruences, refer to \cite[p. 62, Theorem 2.17]{NZM},
$(r-1)i\equiv -t\,(\mod n)$ has a solution in the additive  group $\ZZ_n$ of
integers modulo  $n$ if and only if $\gcd(r-1,n)$ is a divisor of $t$; in this case, the congruence has exactly $\gcd(r-1,n)$  solutions (in $\ZZ_n$). It follows that
\[
|\fix_{V_n}(a_{n,r,t})|=\begin{cases}
\gcd(r-1,n) &\text{ if }\gcd(r-1,n)\mid t\\
0 &\text{ otherwise }.\\
\end{cases}
\]
Using \eqref{eq-a_{r,t}^s-1}, for $s\geqslant 1$, we have
\begin{equation}\label{eq-Vn-art}
 |\fix_{V_n}(a_{n,r,t}^s)|= |\fix_{V_n}(a_{n,r^s,tS_s(r)})|=\begin{cases}
\gcd(r^s-1,n) &\text{ if }\gcd(r^s-1,n)\mid tS_s(r)\\
0 &\text{ otherwise }.\\
\end{cases}
\end{equation}

\begin{lemma}\label{lem-Vn-(n,t)}
Let  $a_{n,r,t}, \, a_{n,r,t'}\in \Aut(\D_{2n})$ with $\gcd(t,n)=\gcd(t',n)$.
Then $\cn(V_n, a_{n,r,t})=\cn(V_n, a_{n,r,t'})$.
 \end{lemma}
 \proof
 By Lemma \ref{lem-kappa-eq-1}(i), $\kappa(n,r,t)=\kappa(n,r,t')$. In addition, it is easy to see that $\gcd(r^s-1,n)\mid tS_s(r)$ if and only if  $\gcd(r^s-1,n)\mid t'S_s(r)$. Then, combining \eqref{eq-cV_n} and \eqref{eq-Vn-art}, the lemma follows.
 \qed

 \begin{lemma}\label{lem-Vn}
Let   $a_{n,r,t}\in \Aut(\D_{2n})$. Assume that $n>1$ is square-free.
Then
\[
 \cn(V_n, a_{n,r,t})=\sum_{
\substack{
  d\mid n\text{ and } \ell\mid \frac{|r|_n}{\gcd\left(\kappa(d,r,t),|r|_n\right)} \\
  \gcd(r^{\ell\kappa(d,r,t)}-1,n)=d
}
}
\frac{d}{\kappa(n,r,t)}\phi\left(\frac{|r|_n}{\ell \gcd\left(\kappa(d,r,t),|r|_n\right)}\right).
\]
\end{lemma}
\proof
In view of \eqref{eq-cV_n} and \eqref{eq-Vn-art}, we first determine those divisors $s$ of $|a_{n,r,t}|=\kappa(n,r,t)$ which satisfy $\gcd(r^s-1,n)\mid tS_s(r)$. Pick such an $s$, and let $d=\gcd(r^s-1,n)$. Then $|r|_d$ is a divisor of $s$, and so $s=k|r|_d$, where $ k\mid \frac{\kappa(n,r,t)}{|r|_d}$. By Lemma~\ref{lem-kappa-eq-1}(iii),
\[tS_{s}(r) \equiv ktS_{|r|_d}(r)\,(\mod d), \]
and so $d$ is a divisor of $ktS_{|r|_d}(r)$. Then $k$ has a divisor
$\frac{d}{\gcd(d, tS_{|r|_d}(r))}$. Write
\[k:=\frac{d\ell}{\gcd(d, tS_{|r|_d}(r))}.\]
Recalling that $n$ is square-free, $\gcd(d,\frac{n}{d})=1$, and so $\frac{\kappa(n,r,t)}{\kappa(d,r,t)}$ is an integer, see Lemma \ref{lem-kappa}.
Since $k$ is a divisor of $\frac{\kappa(n,r,t)}{|r|_d}$, we have
\[ \ell\,\Big |\,\frac{\kappa(n,r,t)}{\kappa(d,r,t)}.\]
Thus,
\[s=\frac{\ell d|r|_d}{\gcd(d, tS_{|r|_d}(r))}=\ell{\kappa(d,r,t)} \mbox{ with } \ell\,\Big |\,\frac{\kappa(n,r,t)}{\kappa(d,r,t)}.\]
 Then we deduce from \eqref{eq-cV_n} and \eqref{eq-Vn-art} that
\[ \begin{split}
\cn(V_n, a_{n,r,t})&{=} \frac{1}{\kappa(n,r,t)}\sum_{
\substack{
  s\mid \kappa(n,r,t) \\
  \gcd(r^s-1,n)\mid tS_{s}(r)
}
}\phi\left( \frac{  \kappa(n,r,t)}{s}\right)  \gcd(r^s-1,n)\\
&=\frac{1}{\kappa(n,r,t)}
\sum_{d\mid n}d
\sum_{
\substack{
  s\mid \kappa(n,r,t) \\
  \gcd(r^s-1,n)=d\\
  d\mid tS_{s}(r)
}
}\phi\left( \frac{  \kappa(n,r,t)}{s}\right) \\
&=\frac{1}{\kappa(n,r,t)}
\sum_{d\mid n}d
\sum_{
\substack{
  \ell\mid \frac{\kappa(n,r,t)}{\kappa(d,r,t)} \\
  \gcd(r^{\ell\kappa(d,r,t)}-1,n)=d\\
  d\mid tS_{\ell\kappa(d,r,t)}(r)
}
}\phi\left( \frac{  \kappa(n,r,t)}{\ell\kappa(d,r,t)}\right) .
\end{split}\]

Assume that $d$ is a divisor of $n$, and   $\ell$ is a divisor of $\frac{\kappa(n,r,t)}{\kappa(d,r,t)}$. Then we obtain from Lemma~\ref{lem-kappa-eq-1}(iii) that
\[ tS_{\ell\kappa(d,r,t)}(r) {\equiv}  \ell tS_{\kappa(d,r,t)}(r) \,(\mod d).\]
Noting that $\kappa(d,r,t)$ is the order of $a_{d,r,t}\in \Aut(\D_{2d})$,
by Lemma \ref{lem-order-a-rt}, \[tS_{\kappa(d,r,t)}(r)\equiv 0\,(\mod d),\] and so
$d\mid tS_{\ell\kappa(d,r,t)}(r)$. In particular, the condition
$d\mid tS_{\ell\kappa(d,r,t)}(r)$ in the summation above is redundant.
Moreover,   $r^{\ell\kappa(d,r,t)}-1$ has a
divisor $r^{|r|_d}-1$. This implies that $d$ is divisor of $r^{\ell\kappa(d,r,t)}-1$, and so $d$  is a common divisor of $n$ and $r^{\ell\kappa(d,r,t)}-1$.
Write $n=dm$.
Recalling that $n$ is square-free, it follows that
  \[
  d=\gcd\left(r^{\ell\kappa(d,r,t)}-1,n\right) \Leftrightarrow \gcd\left(r^{\ell\kappa(d,r,t)}-1,m\right)=1.
  \]
In particular, if $d=\gcd\left(r^{\ell\kappa(d,r,t)}-1,n\right)$ then $\gcd(r-1,m)=1$;
in this case,  $\kappa \left( m,r, t \right)=|a_{m,r, t}|=|r|_m$ by Lemma \ref{lem-order-a-rt}, and so, by Lemma~\ref{lem-kappa},
\[\kappa(n,r,t)=\frac{|r|_n\kappa(d,r,t)\kappa(m,r,t)}{\gcd\left(\kappa(d,r,t),|r|_n\right)
\gcd\left(\kappa(m,r,t),|r|_n\right)}=\frac{|r|_n\kappa(d,r,t)}{\gcd\left(\kappa(d,r,t),|r|_n\right)}.\]

The argument above implies that
\[
\begin{split}
\cn(V_n, a_{n,r,t})&=\frac{1}{\kappa(n,r,t)}
\sum_{d\mid n}d
\sum_{
\substack{
  \ell\mid \frac{\kappa(n,r,t)}{\kappa(d,r,t)} \\
  \gcd(r^{\ell\kappa(d,r,t)}-1,n)=d\\
  d\mid tS_{\ell\kappa(d,r,t)}(r)
}
}\phi\left( \frac{  \kappa(n,r,t)}{\ell\kappa(d,r,t)}\right) \\
&=\sum_{
\substack{
  d\mid n \text{ and } \ell\mid \frac{|r|_n}{\gcd\left(\kappa(d,r,t),|r|_n\right)}  \\
  \gcd(r^{\ell\kappa(d,r,t)}-1,n)=d
}
}
\frac{d}{\kappa(n,r,t)}\phi\left(\frac{|r|_n}{\ell \gcd\left(\kappa(d,r,t),|r|_n\right)}\right).
\end{split}
\]
This completes the proof.
\qed

By \eqref{eq-two-part-c(Dn,art)} and Lemmas~\ref{lem-Un(d)}--\ref{lem-Vn}, the following result is true.

\begin{theorem}
\label{thm-c(D2n, art)}
Let   $a_{n,r,t'}\in \Aut(\D_{2n})$. Assume that $n>1$ is square-free, and let $t=\gcd(t',n)$.
Then
\[ \cn(\D_{2n}^\#, a_{n,r,t'})=\cn(U_n, a_{n,r,0})+\cn(V_n, a_{n,r,t}),\]
where $\cn(U_n, a_{n,r,0})$ and $\cn(V_n, a_{n,r,t})$ are given as in Lemmas~\ref{lem-Un(d)} and Lemma~\ref{lem-Vn}, respectively.
\end{theorem}

\begin{theorem}\label{th-dciD2n}
  Assume that $\D_{2n}$ is a DCI-group with $n>2$. Then the number of non-isomorphic Cayley digraphs
  on $\D_{2n}$ is equal to
  \[\on(2^{\D_{2n}^\#}, \Aut(\D_{2n}))=\frac{1}{n\phi(n)}\sum_{t>0 \text{ and } t\mid n}\phi\left( \frac{n}{t}\right) \sum_{\substack{
                              1\leqslant r\leqslant n-1 \\
                              \gcd(r,n)=1
  }
  }2^{\cn(U_n, a_{n,r,0})+\cn(V_n, a_{n,r,t})},\]
  where $\cn(U_n, a_{n,r,0})$ and $\cn(V_n, a_{n,r,t})$ are given as in Lemmas~\ref{lem-Un(d)} and Lemma~\ref{lem-Vn}, respectively.
\end{theorem}
\proof Since $\D_{2n}$ is a DCI-group  with $n>2$, it follows that $n$ is odd and square-free, see \cite[Corollary 1.2]{XieFK24}. By Lemma \ref{enlemma2-1}, the number of non-isomorphic Cayley digraphs
  on $\D_{2n}$ is equal to $\on(2^{D_{2n}^\#}, \Aut(\D_{2n}))$. Then   the result follows from Theorem \ref{thm-c(D2n, art)}.
\qed

 \vskip 20pt

\section{Cayley digraphs on $\D_{6p}$}\label{sect=D6p}

We continue with the notation in Section \ref{sect=autoD2n}, and let $n=3p$, where $p>3$ is an  odd prime. In the following, we let  $d,\,t\in \{1,3,p,3p\}$, and $1\leqslant r\leqslant 3p-1$ with $\gcd(r,n)=1$.

By Lemma \ref{lem-Un(d)},
\begin{equation*}
\cn(U_{3p}, a_{3p,r,0})=\sum_{d>1 \text{ and }d \mid 3p} \frac{\phi(d)}{|r|_d}=\frac{2}{|r|_3}+\frac{p-1}{|r|_p}+\frac{2(p-1)}{|r|_{3p}}.
\end{equation*}
By Lemma \ref{r_n-d}, $|r|_{3p}=\frac{|r|_3|r|_p}{\gcd(|r|_3,|r|_p)}$. Noting that $r\equiv 1,\,2\,(\mod 3)$, we have
\begin{equation}\label{r-3p}
 |r|_{3p}=
 \begin{cases}
   |r|_p & \mbox{if }  r\equiv 1\,(\mod 3)\\
   \frac{2|r|_p}{\gcd(2,|r|_p)}  & \mbox{if }  r\equiv 2\,(\mod 3).
 \end{cases}
\end{equation}
Then $\cn(U_{3p}, a_{3p,r,0})$ is known as follows.
\begin{lemma}\label{eq-c(U3p,art)}
Let $t\in \{1,3,p,3p\}$, and $1\leqslant r\leqslant 3p-1$ with $\gcd(r,n)=1$. Then
\[
\cn(U_{3p}, a_{3p,r,0})=
 \begin{cases}
   2+\frac{3(p-1)}{|r|_p} & \mbox{if }  r\equiv 1\,(\mod 3)\\
   1+\frac{2(p-1)}{|r|_p}  & \mbox{if }  r\equiv 2\,(\mod 3)\mbox{ and $|r|_p$ is odd}
   \\
   1+\frac{3(p-1)}{|r|_p}  & \mbox{if }  r\equiv 2\,(\mod 3)\mbox{ and $|r|_p$ is even}.
 \end{cases}
\]
\end{lemma}

Put \[\Lambda(d,r,t):=\{\ell>0\mid \ell\text{ divides } \frac{|r|_{3p}}{\gcd\left(\kappa(d,r,t),|r|_{3p}\right)} \text{ and } \gcd(r^{\ell\kappa(d,r,t)}-1,3p)=d\}.\]
For a positive integer $m$, put
\[\Delta(m):=\{0<\ell<m\mid \ell \mbox{ divides } m\}.\]
Recalling the definition of $\kappa(n,r,t)$, we see that $\kappa(1,r,t)=1$, and so $\frac{|r|_{3p}}{\gcd\left(\kappa(1,r,t),|r|_{3p}\right)}=|r|_{3p}$. Then it follows from \eqref{eq-c(U3p,art)} that
\begin{equation}\label{quotient-1}
 \frac{|r|_{3p}}{\gcd\left(\kappa(1,r,t),|r|_{3p}\right)}=
 \begin{cases}
   |r|_p & \mbox{if }  r\equiv 1\,(\mod 3)\\
   \frac{2|r|_p}{\gcd(2,|r|_p)}  & \mbox{if }  r\equiv 2\,(\mod 3).
 \end{cases}
\end{equation}
Next we derive the formula of $\Lambda(1,r,t)$.

\begin{lemma}\label{V3p-1-t}
Let $t\in \{1,3,p,3p\}$, and $1\leqslant r\leqslant 3p-1$ with $\gcd(r,n)=1$. Then \[\Lambda(1,r,t)
=\begin{cases}
   \emptyset & \mbox{if } r\equiv  1\,(\mod 3)\\
   \{\ell\in \Delta(|r|_p)\mid   \ell \mbox{ is odd}\}
   & \mbox{if } r\equiv  2\,(\mod 3).
 \end{cases}
\]
\end{lemma}
\proof It is clear that \[\Lambda(1,r,t)=\{\ell>0\mid \ell\mbox{ divides } |r|_{3p} \text{ and } \gcd(r^{\ell}-1,3p)=1\}.\]
Let $\ell$ divide $|r|_{3p}$. Then $\ell\in \Lambda(1,r,t)$ if and only if
 $\gcd(r^{\ell}-1,3p)=1$ or, equivalently,
  neither $|r|_3$ nor $|r|_p$ is a divisor of $\ell$.
In particular, if $|r|_3=1$, i.e., $r\equiv  1\,(\mod 3)$, then $\Lambda(1,r,t)=\emptyset$.
Suppose that $r\equiv  2\,(\mod 3)$, and so  $|r|_3=2$. Then  $\ell\in \Lambda(1,r,t)$ if and only if $\ell$ is odd and indivisible by $|r|_p$.
By \eqref{r-3p},  $|r|_{3p}=\frac{2|r|_p}{\gcd(2,|r|_p)}$.
It follows that  $\ell$ is an odd divisor of $|r|_{3p}$ if and only if $\ell$ is an odd divisor of $|r|_{p}$. Then the lemma follows.
\qed

\begin{lemma}\label{V3p-3-t}
Let $t\in \{1,3,p,3p\}$, and $1\leqslant r\leqslant 3p-1$ with $\gcd(r,n)=1$. Then \[\Lambda(3,r,t)
=\begin{cases}
\Delta\left(|r|_p\right)  & \mbox{if } r\equiv  1\,(\mod 3)\mbox{ and }t\in \{3,3p\}\\
   \Delta\left(\frac{|r|_p}{\gcd(3,|r|_p)}\right)  & \mbox{if } r\equiv  1\,(\mod 3)\mbox{ and }t\in \{1,p\}\\
   \Delta\left(\frac{|r|_p}{\gcd(2,|r|_p)}\right)   & \mbox{if } r\equiv  2\,(\mod 3).
 \end{cases}
\]
\end{lemma}
\proof
Calculation shows that
\[\kappa(3,r,t)=
\begin{cases}
  1 & \mbox{if } t\in \{3,3p\}\mbox{ and }  r\equiv  1\,(\mod 3) \\
  2 & \mbox{if } t\in \{3,3p\}\mbox{ and }  r\equiv  2\,(\mod 3) \\
  3 & \mbox{if } t\in \{1,p\} \mbox{ and }  r\equiv  1\,(\mod 3)\\
  2 & \mbox{if } t\in \{1,p\} \mbox{ and }  r\equiv  2\,(\mod 3) .
\end{cases}
\]

Assume that $ r\equiv  1\,(\mod 3)$.
Then $ \gcd(r^{\ell\kappa(3,r,t)}-1,3p)=3$ if and only if $|r|_p$ is not a divisor of $\ell\kappa(3,r,t)$.
 By \eqref{r-3p},  $|r|_{3p}= |r|_p$, we have
\begin{equation}\label{quotient-2-1}
  \frac{|r|_{3p}}{\gcd\left(\kappa(3,r,t),|r|_{3p}\right)}=
\begin{cases}
  |r|_p
  & \mbox{if } t\in \{3,3p\}, \\
   \frac{|r|_p}{\gcd(3,|r|_p)}
  & \mbox{if } t\in \{1,p\}.
\end{cases}
\end{equation}
Then
\begin{equation}\label{3-l-1}
\ell\in \Lambda(3,r,t)\Leftrightarrow
\begin{cases}
 \ell\in \Delta(|r|_p)
  & \mbox{if } t\in \{3,3p\}, \\
  \ell\in \Delta\left(\frac{|r|_p}{\gcd(3,|r|_p)}\right)
  & \mbox{if } t\in \{1,p\}.
\end{cases}
\end{equation}

Assume that $ r\equiv  2\,(\mod 3)$. Then $\kappa(3,r,t)=2=|r|_3$. Then
$ \gcd(r^{\ell\kappa(3,r,t)}-1,3p)=3$ if and only if
$\ell\kappa(3,r,t)$ is  indivisible by $|r|_p$. By \eqref{r-3p},
 $|r|_{3p}= \frac{2|r|_p}{\gcd(2,|r|_p)}$, we have

 \begin{equation}\label{quotient-2-2}
  \frac{|r|_{3p}}{\gcd\left(\kappa(3,r,t),|r|_{3p}\right)}=\frac{|r|_p}{\gcd(2,|r|_p)}.
 \end{equation}
Then
\begin{equation}\label{3-l-2}\ell\in \Lambda(3,r,t)\Leftrightarrow
  \ell\in \Delta\left(\frac{|r|_p}{\gcd(2,|r|_p)}\right).
\end{equation}

Finally, the lemma follows from \eqref{3-l-1} and \eqref{3-l-2}.
\qed

\begin{lemma}\label{V3p-p-t}
Let $t\in \{1,3,p,3p\}$, and $1\leqslant r\leqslant 3p-1$ with $\gcd(r,n)=1$. Then \[\Lambda(p,r,t)
=\begin{cases}
  \{1\}  & \mbox{if } r\equiv  2\,(\mod 3) \mbox{ and } |r|_p\mbox{ is odd}\\
  \emptyset   & \mbox{otherwise }.
 \end{cases}
\]
\end{lemma}
\proof
Calculation shows that
\[\kappa(p,r,t)=
\begin{cases}
  p  & \mbox{if } t\in \{1,3\}\mbox{ and }  r\equiv  1\,(\mod p) \\
  |r|_p & \mbox{otherwise}.
\end{cases}
\]
Recalling that $|r|_{3p}=\frac{|r|_3|r|_p}{\gcd(|r|_3,|r|_p)}$, we have
\begin{equation}\label{quotient-3}
  \frac{|r|_{3p}}{\gcd\left(\kappa(p,r,t),|r|_{3p}\right)}=
\begin{cases}
  |r|_3
  & \mbox{if } t\in \{1,3\}\mbox{ and }  r\equiv  1\,(\mod p)\\
  \frac{|r|_3}{\gcd(|r|_3,|r|_p)}
  &  \mbox{otherwise}.
\end{cases}
\end{equation}
Note that $ \gcd(r^{\ell\kappa(p,r,t)}-1,3p)=p$ if and only if   $\ell\kappa(p,r,t)$ is divisible by $|r|_p$ but not by $|r|_3$. In particular, if $\ell\in \Lambda(p,r,t)$ then
 $|r|_3=2$ and $\ell\kappa(p,r,t)$ is odd, and so $\ell=1$ and $\kappa(p,r,t)$ is odd.
 Then the lemma follows.
 \qed

 Recall that $|r|_{3p}$ is a divisor of $\kappa(3p,r,t)$, see Lemma \ref{lem-order-a-rt}.
We have the following lemma.

\begin{lemma}\label{V3p-3p-t}
Let $t\in \{1,3,p,3p\}$, and $1\leqslant r\leqslant 3p-1$ with $\gcd(r,n)=1$. Then $\Lambda(3p,r,t)=\{1\}$.
\end{lemma}

Put
\[\mathcal{S}(d,r,t):=\sum_{\ell\in \Lambda(d,r,t)
}
d\phi\left(\frac{|r|_{3p}}{\ell \gcd\left(\kappa(d,r,t),|r|_{3p}\right)}\right).\]
By Lemma \ref{V3p-3p-t}, we have the following lemma.
\begin{lemma}\label{S-3prt}
Let $t\in \{1,3,p,3p\}$, and $1\leqslant r\leqslant 3p-1$ with $\gcd(r,n)=1$. Then
$\mathcal{S}(3p,r,t)=3p$.
\end{lemma}

By Lemma \ref{V3p-p-t} and \eqref{quotient-3}, the following lemma holds.
\begin{lemma}\label{S-prt}
Let $t\in \{1,3,p,3p\}$, and $1\leqslant r\leqslant 3p-1$ with $\gcd(r,n)=1$. Then \[\mathcal{S}(p,r,t)
=\begin{cases}
 p  & \mbox{if } r\equiv  2\,(\mod 3) \mbox{ and } |r|_p\mbox{ is odd}\\
  0   & \mbox{otherwise}.
 \end{cases}
\]
\end{lemma}

\begin{lemma}\label{S-1rt}
Let $t\in \{1,3,p,3p\}$, and $1\leqslant r\leqslant 3p-1$ with $\gcd(r,n)=1$. Then \[\mathcal{S}(1,r,t)
=\begin{cases}
  0  & \mbox{if } r\equiv  1\,(\mod 3)\\
  |r|_p-1   &  \mbox{if } r\equiv  2\,(\mod 3)\mbox{ and $|r|_p$ is odd}\\
 \frac{|r|_p}{2}   &  \mbox{if } r\equiv  2\,(\mod 3)\mbox{ and $|r|_p$ is even}.
 \end{cases}
\]
\end{lemma}
\proof
By Lemma~\ref{V3p-1-t} and \eqref{quotient-1},
\[
\begin{split}
\mathcal{S}(1,r,t)&=\sum_{\ell\in \Lambda(1,r,t)
}
\phi\left(\frac{|r|_{3p}}{\ell \gcd\left(\kappa(1,r,t),|r|_{3p}\right)}\right)\\
&=\begin{cases}
   0  & \mbox{if } r\equiv  1\,(\mod 3)\\
  \sum_{\mbox{\tiny odd }\ell\in \Delta(|r|_p)
}\phi\left(\frac{2|r|_{p}}{\ell\gcd(2,|r|_{p})}\right) & \mbox{if } r\equiv  2\,(\mod 3).
 \end{cases}
\end{split}
\]
Suppose that $r\equiv  2\,(\mod 3)$. If $|r|_p$ is odd
then
\begin{align*}
\mathcal{S}(1,r,t)
&= \sum_{\mbox{\tiny odd }\ell\in \Delta(|r|_p)
}\phi\left(\frac{2|r|_{p}}{\ell\gcd(2,|r|_{p})}\right)\\
&= \sum_{\mbox{\tiny odd }\ell\in \Delta(|r|_p)
}\phi\left(\frac{2|r|_{p}}{\ell}\right)
=\sum_{\mbox{\tiny odd }\ell\in \Delta(|r|_p)
}\phi\left(\frac{|r|_{p}}{\ell}\right)=|r|_p-1.
\end{align*}
If $|r|_p$ is even
then
\[\mathcal{S}(1,r,t)= \sum_{\mbox{\tiny odd }\ell\in \Delta(|r|_p)
}\phi\left(\frac{2|r|_{p}}{\ell\gcd(2,|r|_{p})}\right)= \sum_{\mbox{\tiny odd }\ell\in \Delta(|r|_p)
}\phi\left(\frac{|r|_{p}}{\ell}\right)=\frac{|r|_p}{2}.\]
Thus, the lemma is proved.
\qed

\begin{lemma}\label{S-3rt}
Let $t\in \{1,3,p,3p\}$, and $1\leqslant r\leqslant 3p-1$ with $\gcd(r,n)=1$. Then \[\mathcal{S}(3,r,t)
=\begin{cases}
  3(|r|_p-1)  & \mbox{if } r\equiv  1\,(\mod 3) \mbox{ and }t\in\{3,3p\}\\
  3(|r|_p-1)   &  \mbox{if } r\equiv  1\,(\mod 3),\,t\in\{1,p\} \mbox{ and } 3\nmid |r|_p
  \\
 |r|_p-3  &  \mbox{if } r\equiv  1\,(\mod 3),\,t\in\{1,p\} \mbox{ and }3\mid |r|_p\\
 3(|r|_p-1)  &  \mbox{if } r\equiv  2\,(\mod 3)\mbox{ and $|r|_p$ is odd}
\\
3\left(\frac{|r|_p}{2}-1\right)  &  \mbox{if } r\equiv  2\,(\mod 3)\mbox{ and $|r|_p$ is even}.
 \end{cases}
\]
\end{lemma}
\proof
If $r\equiv  1\,(\mod 3)$, then we deduce from \eqref{quotient-2-1} and \eqref{3-l-1} that
 \[
 \begin{split}
 \mathcal{S}(3,r,t)
&=\sum_{\ell\in \Lambda(3,r,t)
}
3\phi\left(\frac{|r|_{3p}}{\ell \gcd\left(\kappa(3,r,t),|r|_{3p}\right)}\right)\\
&=\begin{cases}
 \sum_{\ell\in \Delta(|r|_p)}3\phi(\frac{|r|_p}{\ell}) =3(|r|_p-1)
  & \mbox{if } t\in \{3,3p\} \\
   \sum_{\ell\in \Delta\left(\frac{|r|_p}{\gcd(3,|r|_p)}\right)}3\phi\left(\frac{|r|_p}{\ell\gcd(3,|r|_p)}\right)
   =3\left(\frac{|r|_p}{\gcd(3,|r|_p)}-1\right)
  & \mbox{if } t\in \{1,p\}.
\end{cases}
\end{split}
\]
If $r\equiv  2\,(\mod 3)$, then we obtain from \eqref{quotient-2-2} and \eqref{3-l-2} that
 \[
 \begin{split}
 \mathcal{S}(3,r,t)
&=\sum_{\ell\in \Lambda(3,r,t)
}
3\phi\left(\frac{|r|_{3p}}{\ell \gcd\left(\kappa(3,r,t),|r|_{3p}\right)}\right)\\
&=\sum_{\ell\in \Delta\left(\frac{|r|_p}{\gcd(2,|r|_p)}\right)}3\phi\left(\frac{|r|_p}{\ell\gcd(2,|r|_p)}\right)
   =3\left(\frac{|r|_p}{\gcd(2,|r|_p)}-1\right).
\end{split}
\]
Thus, the lemma follows.
\qed

\begin{lemma}\label{eq-c(V3p,art)}
Let $t\in \{1,3,p,3p\}$, and $1\leqslant r\leqslant 3p-1$ with $\gcd(r,n)=1$. Then
 \[
\cn(V_{3p}, a_{3p,r,t})=\begin{cases}
\gcd(3p,t) & \mbox{if } r=1,\\
  3+\frac{3(p-1)}{|r|_p}  & \mbox{if } 1\ne r\equiv  1\,(\mod 3) \mbox{ and }t\in\{3,3p\}\\
 1+\frac{p-1}{|r|_p} &  \mbox{if } 1\ne r\equiv  1\,(\mod 3),\,t\in\{1,p\} \mbox{ and } 3\nmid |r|_p
  \\
   1+\frac{3(p-1)}{|r|_p} &  \mbox{if } 1\ne r\equiv  1\,(\mod 3),\,t\in\{1,p\} \mbox{ and } 3\mid |r|_p
  \\
 2\gcd(p,t) &  \mbox{if } r\equiv  2\,(\mod 3),\,r\equiv  1\,(\mod p)\\
  2+\frac{2(p-1)}{|r|_p}  &  \mbox{if } r\equiv  2\,(\mod 3),\,r\not\equiv  1\,(\mod p)\mbox{ and $|r|_p$ is odd}
\\
2+\frac{3(p-1)}{|r|_p}   &  \mbox{if } r\equiv  2\,(\mod 3),\,r\not\equiv  1\,(\mod p)\mbox{ and $|r|_p$ is even}.
 \end{cases}
\]
 \end{lemma}
 \proof
 By \eqref{r-3p},
 \[|r|_{3p}=
 \begin{cases}
   |r|_p & \mbox{if }  r\equiv 1\,(\mod 3)\\
    2|r|_p  & \mbox{if }  r\equiv 2\,(\mod 3)\mbox{ and $|r|_p$ is odd}\\
    |r|_p  & \mbox{if }  r\equiv 2\,(\mod 3)\mbox{ and $|r|_p$ is even}.
 \end{cases}
 \]
 Recalling that $\kappa(3p,r,t) =\frac{3p|r|_{3p}}{\gcd\left(3p, tS_{|r|_{3p}}(r)\right)}$, we have
 \[
  \kappa(3p,r,t) =
 \begin{cases}
 \frac{3p}{\gcd(3p,t)} & \mbox{if }  r=1\\
  \frac{3|r|_p}{\gcd(3,t|r|_p)} & \mbox{if }  1\ne r\equiv 1\,(\mod 3)\\
  \frac{2p}{\gcd(p,t)} & \mbox{if }  r\equiv 2\,(\mod 3)\mbox{ and }  r\equiv 1\,(\mod p)\\
   \frac{2|r|_p}{\gcd(2,|r|_p)}  & \mbox{if }  r\equiv 2\,(\mod 3)\mbox{ and  }  r\not\equiv 1\,(\mod p).
 \end{cases}
 \]
By Lemma \ref{lem-Vn},
\[
\begin{split}
\cn(V_{3p}, a_{3p,r,t})&=\frac{1}{\kappa(3p,r,t)}\sum_{d\in\{1,3,p,3p\}}\sum_{\ell\in \Lambda(d,r,t)
}
d\phi\left(\frac{|r|_{3p}}{\ell \gcd\left(\kappa(d,r,t),|r|_{3p}\right)}\right)\\
&=\frac{1}{\kappa(3p,r,t)}\left(\mathcal{S}(3p,r,t)+\mathcal{S}(p,r,t)+\mathcal{S}(3,r,t)
+\mathcal{S}(1,r,t)\right).
\end{split}
\]
Then, using Lemmas \ref{S-3prt}-\ref{S-3rt}, the lemma follows  from a straightforward calculation.
 \qed

 By \cite[Theorem 1]{DobsonMS2015}, $\D_{6p}$ is a DCI-group if and only if $p \geqslant 5$.
Then the number of non-isomorphic Cayley digraphs on $\D_{6p}$ (with $p>3$) is given in the following result, which follows  from Theorem \ref{th-dciD2n}, Lemmas \ref{eq-c(U3p,art)} and \ref{eq-c(V3p,art)}.

\begin{theorem}\label{th-dciD6p}
 Let $p>3$ be a prime.  Then the number of non-isomorphic Cayley digraphs
  on $\D_{6p}$ is equal to
  \[
  \begin{split}
  \on(2^{\D_{6p}^\#}, \Aut(\D_{6p}))&=\frac{2^{4p-2}}{3p(p-1)}\left(2^{2p}+5\right)+\frac{2^{2p}}{p}\left(2^{p}+1\right)\\
  &+\frac{8}{3(p-1)}\sum_{ |r|_3=1\text{ and }3\nmid |r|_p}\left(2^{\frac{4(p-1)}{|r|_p}}+2^{1+\frac{6(p-1)}{|r|_p}}\right)+\frac{8}{p-1}\sum_{|r|_3=1\text{ and }3\mid |r|_p}2^{\frac{6(p-1)}{|r|_p}}\\
  &+\frac{4}{p-1}\left(\sum_{|r|_3=2\text{ and }2\nmid |r|_p\text{ and }|r|_p>1}2^{\frac{4(p-1)}{|r|_p}}+\sum_{|r|_3=2\text{ and }2\mid |r|_p}2^{\frac{6(p-1)}{|r|_p}}\right),
  \end{split}
   \]
   where $r$ runs over $\{r\mid 1<r<3p, \gcd(r,3p)=1\}$.
\end{theorem}

\end{document}